\documentclass[11pt,a4paper]{amsart}
\usepackage[english]{babel}
\usepackage{amssymb,xypic,amscd}
\CompileMatrices

\newtheorem{defn}{Definition}
\newtheorem{thm}[defn]{Theorem}

\newtheorem{cor}[defn]{Corollary}
\newtheorem{lem}[defn]{Lemma}
\newtheorem{prop}[defn]{Proposition}

\theoremstyle{remark}
\newtheorem{rem}[defn]{Remark}
\theoremstyle{remark}
\newtheorem{exam}{Example}

\numberwithin{equation}{section} \numberwithin{defn}{section}

\renewcommand\sp{\operatorname{Spec}}

\newcommand\Ext{\operatorname{Ext}}
\renewcommand\hom{\operatorname{Hom}}
\newcommand\Ker{\operatorname{Ker}}

\newcommand\Pic{\operatorname{Pic}}
\newcommand\Apl{\operatorname{Map}}
\newcommand\Coker{\operatorname{Coker}}

\newcommand\res{\operatorname{Res}}
\renewcommand\deg{\operatorname{deg}}

\newcommand\limpl[1]{\underset{#1}\varprojlim\,}

\renewcommand\o{{\mathcal O}}
\renewcommand\L{{\mathcal L}}
\newcommand\B{{\mathcal B}}

\newcommand\M{{\mathcal M}}

\renewcommand\tilde{\widetilde}

\renewcommand\lim{\limpl{A\in\B}}
\newcommand\beq{
      \setcounter{equation}{\value{defn}}\addtocounter{defn}1
      \begin{equation}}
\def\mod #1/#2{\kern.06em{\raise1.2pt\hbox{$#1$}}/
         {\raise-1.2pt\hbox{$#2$}}}

\begin{document}

\title[Generalized Local Symbols]{GEOMETRIC CONSTRUCTION OF GENERALIZED \\ LOCAL SYMBOLS ON ALGEBRAIC CURVES}
\author[Fernando Pablos Romo]
{Fernando Pablos Romo \\
Departamento de Matem\'aticas \\ Universidad de Salamanca}

\address{Departamento de
Matem\'aticas, Universidad de Salamanca,  Plaza de la Merced 1-4\\
Salamanca 37008. Spain.}

\thanks{
This work is partially supported by the DGESYC research contract
n. BFM2003-00078 and Castilla y Le\'on Regional Government
 contract SA071/04.}

\email{fpablos@usal.es}

\begin{abstract}    The aim of this work is to provide a construction of generalized
local symbols on algebraic curves as morphisms of group schemes.
From a closed point of a complete, irreducible and non-singular
curve $C$ over a perfect field $k$ as the only data, using theta
groups over Picard schemes of curves, we offer a geometric
construction that allows us to define generalizations of the tame
symbol and the Hilbert norm residue symbol.
\end{abstract}

\maketitle


\setcounter{tocdepth}1


\section{Introduction}

    In 1959 J.-P. Serre \cite{Se} introduced local symbols on curves in order
to study the factorization of a rational morphism of complete
curves through a generalized Jacobian of Rosenlicht (\cite{Ro1},
\cite{Ro2}). The additive local symbol (the classic residue on
curves) and the multiplicative local symbol are examples of these
symbols.

    Given a complete, irreducible and non-singular curve $C$ over an algebraically
closed field $k$ and a closed point $p\in C$, for each $f,g\in
{\Sigma}_C^*$ (the multiplicative group  of the function field
${\Sigma}_C$ of $C$), the expression of the multiplicative local
symbol is:
$$(f,g)_p = (-1)^{v_p(f)\cdot v_p(g)} \frac
{f^{v_p(g)}}{g^{v_p(f)}} (p) \in k^*\, .$$

    A few years later, in 1971, J. Milnor \cite{Mi} defined a generalization of
the multiplicative local symbol, the tame symbol $d_v$, associated
with a discrete valuation $v$ on a field $F$. Explicitly, if $A_v$
is the valuation ring,
 $p_v$ is the unique maximal ideal and $k_v = \mod {A_v}/{p_v}$ is
 the residue class field, Milnor defined $d_v \colon F^*\times F^* \to k_v^*$
by $$d_v(x,y) = (-1)^{v(x)\cdot v(y)} \frac {x^{v(y)}}{y^{v(x)}}
(\text {mod }p_v).$$

        Moreover, when $C$ is defined over a finite field that contains the
$m^{th}$-roots of unity, the formula of the Hilbert norm residue
symbol is:
$$(f,g)_p =  \left (N_{k(p)/k}
[(-1)^{v_p(f)\cdot v_p(g)} \frac {f^{v_p(g)}}{g^{v_p(f)}}
(p)]\right )^{\frac {q-1}{m}} \in {\mu}_m\, ,$$ \noindent where
$\# k = q$. This expression was offered by H.L. Schmid \cite{Sc},
and we can obtain it from the multiplicative local symbol by using
the norm morphism.

    In recent last years both symbols, the multiplicative local symbol
and the Hilbert norm residue symbol, have been studied by several
authors from different points or view (\cite {ACK}, \cite{Be},
\cite{Co}, \cite{De}, \cite{Ho}, \cite{Pa0} or \cite{Pa}).

    The goal of this paper is to give an algebraic construction
that will allow us to define generalized local symbols as
morphisms of schemes, and by using a closed point of a complete,
irreducible and non-singular curve $C$ over a perfect field $k$,
as the unique datum to make the construction.

    As results, we offer morphisms of group schemes that generalize
the multiplicative local symbol and several symbols that we can
obtain from it, among them the Hilbert norm residue symbol.

    We should remark that generalizations of the tame symbol have been made
by C. Contou-Carr\`{e}re \cite{Co} as a morphism of functors and by
the author as a morphism of schemes (\cite{Pa3}). This paper gives
the relation between the classic datum used to define a local
symbol (a closed point of a curve) and the formal scheme that
appears in \cite{Pa3}. As far as we know, generalized local
symbols, different from the symbols referred to previously, have
not been stated explicitly in the literature. Following the method
of \cite{AP}, we also offer a reciprocity law for these
generalized symbols when we consider $\sp B$-valued points, $B$
being an artinian local finite $k$-algebra.

    The organization of the paper is as follows:

    In Section 2 definitions and basic results are introduced. Thus, Witt's schemes are
defined from a ring scheme, the Jacobian of a cuspidal curve
$C^{n)}$, constructed from a smooth curve $C$, is characterized,
and the group scheme $K[{\mathcal L}^{n)}]$, associated with an
invertible sheaf on the Jacobian $J_{C^{n)}}$, are also defined.

    Section 3 contains the properties of the group schemes
$K[{\tilde {\mathcal L}}]$ over Picard schemes of curves, which
are basic to the geometric construction offered. Group schemes
associated with invertible sheaves over $\Pic C$ and $\Pic C^{n)}$
are defined, and a detailed study of the group $K[({\tilde
{{\mathcal L}_{\Theta}}})^{n)}]$, associated with the invertible
sheaf $({\tilde {{\mathcal L}_{\Theta}}})^{n)}$ over $\Pic C^{n)}$
and constructed from the principal polarization of $J_C$, is made.

    Finally, Section 4 is devoted to constructing, from the scheme
$K[({\tilde {{\mathcal L}_{\Theta}}})^{n)}]$, a formal scheme that
allows us to define morphisms of schemes that generalize the local
symbols referred to above. Similar to \cite{Pa3}, we provide two
generalizations of the multiplicative local symbol: one is
associated with a rational point and the other, which is more
general, is associated with an arbitrary closed point. The first
generalization coincides, when we consider rational points, with
the expression of J.P. Serre and, if we consider $S$-valued points
with $S$ a connected $k$-scheme, with the formula of C.
Contou-Carr\`{e}re. The second generalization offered allows us to
provide new morphisms of schemes, from morphisms of groups schemes
${\mathbb G}_m \to G$, and, in particular, from the characters of
${\mathbb G}_m$. Hence, the Hilbert norm residue symbol is
generalized and the classic formula is again recovered by
considering rational points. To conclude, we offer reciprocity
laws for these generalized symbols when one considers $\sp
B$-valued points, $B$ being an artinian local finite $k$-algebra.

\section{Preliminaries}

        Firstly, we shall recall the concept of $S$-valued points in
the theory of schemes: if $Y$ and $S$ are $k$-schemes, an
$S$-valued point of $Y$ is a morphism of schemes from $S$ to $Y$.
The set of all such is denoted $$Y^{\bullet} (S) =
\hom_{k-sch.}(S,Y)\, .$$

\subsection{Witt's Schemes}

 Let $\mathcal A$ be a ring $k$-scheme: i.e., it is a scheme
 whose functor of points takes values in the category of commutative
 rings with unit element.

       Let us assume that $\mathcal A$ is unipotent a group scheme
 under addition (there exists an isomorphism of group schemes between $\mathcal A$
and a unipotent group scheme). We shall now define the Witt
schemes associated with $\mathcal A$.

       Recall from (\cite{Se}, page 168) that given two
commutative group schemes, $A$ and $B$, $H^2 (A,B)$ is the group
of classes of $2$-cocycles (mod. $2$-coboundaries), and that
$H^2_{reg} (A,B)_s$ is the subgroup of $H^2 (A,B)$ consisting of
classes of regular and symmetric 2-cocycles. It is also known that
each regular and symmetric 2-cocycle determines an extension
between unipotent groups, i.e. $H^2_{reg} (A,B)_s \subseteq \Ext
(A,B)$. Moreover, if the ground field is algebraically closed,
each extension of unipotent groups is induced by a regular
2-cocycle (\cite{Se}, page 171).

        Let us now consider the regular morphism of schemes $f_2\colon {\mathcal A} \times {\mathcal A} \to
 {\mathcal A}$, which defines the multiplication in $\mathcal A$.
 One has that
$f_2$ satisfies the condition of a 2-cocycle: $$f_2(y,z) -
f_2(x+y,z) + f_2(x,y+z) - f_2(x,y) = 0\, ,
$$\noindent where $x,y,z\in {\mathcal A}^{\bullet}(S)$, $S$ being an arbitrary $k$-scheme.

        Hence, $f_2\in H^2_{reg} ({\mathcal A},{\mathcal A})_s$ and it induces an extension
$$0\to {\mathcal A} \to W^2_{\mathcal A}\to {\mathcal A} \to 0\,
$$\noindent such that $W^2_{\mathcal A}\simeq {\mathcal A} \times {\mathcal
A}$ as $k$-schemes, and its group structure is defined by the map
$$m[(x,y),(\bar x,\bar y)] = (x + \bar x,y + \bar y + f_2(x,\bar
x))\, ,$$ for all $x, \bar x, y,\bar y \in {\mathcal
A}^{\bullet}(S)$, with $S$ a $k$-scheme.

        Since $W^2_{\mathcal A}$ is a unipotent group, by induction we obtain the extensions
 $0\to {\mathcal A} \to
W^h_{\mathcal A}\to W^{h-1}_{\mathcal A}\to 0$, determined by the
morphisms of schemes $$f_h\colon W^{h-1}_{\mathcal A}\times
W^{h-1}_{\mathcal A}\to {\mathcal A}\, ,$$\noindent where $$f_h
(x,y) = \sum_{\begin{aligned} i+j &= h \\ i,j&\geq 1
\end{aligned}} x_i\cdot y_j$$\noindent for all $x=(x_i),y=(y_i) \in
(W^{h-1}_{\mathcal A})^{\bullet}(S)$.

        If $z=(z_i)\in (W^{h-1}_{\mathcal
A})^{\bullet}(S)$, it follows from the equality $$f_h(y,z) -
f_h(x+y,z) + f_h(x,y+z) - f_h(x,y) = 0$$ that $f_h$ is a 2-cocycle
and that it determines a factor system that defines the group
structure of $W^{h}_{\mathcal A}$.

\begin{exam}\label{e:witt} Let $K$ be a finite extension of the field $k$.

    From the multiplicative law of $K$,
$m\colon K\times K \to K$,  we have an induced map $K^{\vee}\to
K^{\vee}\oplus K^{\vee}$, where $K^{\vee} = {\hom}_k (K,k)$. Thus,
by computing the symmetric algebra, we get a morphism of
$k$-algebras $$S^{\cdot} (K^{\vee}) \longrightarrow S^{\cdot} (K^{\vee})
\underset {k} \otimes S^{\cdot} (K^{\vee})\, ,$$\noindent which induces
a map $\sp S^{\cdot} (K^{\vee}) \times \sp S^{\cdot} (K^{\vee})
\longrightarrow \sp S^{\cdot} (K^{\vee})$, which, together with the
morphism of schemes induced over $\sp S^{\cdot} (K^{\vee})$ by the
addition map of $K$, determines the structure of a ring scheme on
$\sp S^{\cdot} (K^{\vee})$.

       Furthermore, if $d = \text{ dim}_k K$, there exists an isomorphism
of additive groups $K\simeq k^{d}$, and therefore the ring scheme
$\sp S^{\cdot} (K^{\vee})$ is isomorphic as a group scheme to ${\mathbb
G}_a\times \overset {d}
 \dots \times {\mathbb G}_a$. However, both schemes are not isomorphic as ring schemes.

        Let us now consider a $k$-algebra $A$, and let us set
 $W_A^{n}$ to denote the group $A^n$ together with the operation $(a_1,\dots,a_n) +
(b_1,\dots,b_n) = (c_1,\dots,c_n)$, where $c_i = a_i + b_i +
\sum_{h+k=i} a_h\cdot b_k$.

      Hence, we can now define the Witt scheme $W^h_{\sp
S^{\cdot} (K^{\vee})}$, whose functor of points is $$ [W^h_{\sp S^{\cdot}
(K^{\vee})}]^{\bullet} (S) =W^h_{\Gamma (S,{\o}_S \underset {k}
\otimes K)}\, ,$$\noindent for each $k$-scheme $S$.
\end{exam}

    Moreover, if $W(A)^n_{+}$ is the $n^{th}$-group defined in (\cite{Haz}), one
has that $W(A)^n_{+}$ and $W_A^{n}$ are isomorphic through the map
$$\aligned W(A)^n_{+}
&\longrightarrow  W_A^{n} \\ ({\bar a}_i) &\longrightarrow (b_j)\,
,
\endaligned$$ with $$b_j =\sum \sb{ \overset {i_1<i_2<\dots<i_k}
{i_1 + \dots + i_k = j}} (-1)^k {\bar a}_{i_1}\dots {\bar
a}_{i_k}.$$

    Furthermore, both groups are also isomorphic to
 ${\Gamma}_+^n(A)$, which is the multiplicative subgroup of
$(A[t]/{t^{n+1}})^{*}$ consisting of elements of the form $$\{1 +
a_1t + \dots + a_nt^n\}\, .$$

\subsection{Cuspidal Curve associated with a Smooth Curve}

   Let us now consider a complete, irreducible and non-singular curve $C$ over a perfect field $k$ and with a rational
point. If $p\in C$ is a closed point, ${\mathcal O}_{p}$ is its
local ring, and $k(p)$ is the residue field at $p$, for each
positive integer $n$ one has that there exists the scheme
$$C^{n)} =   \underset {\sp (\mod
{{\o}_p}/{m_p^n})\qquad} {C \coprod \sp k(p)}\, ,$$\noindent such
that:
\begin{itemize}
\item The topological spaces of $C$ and $C^{n)}$ are
homeomorphic.
\item If $\{U_i\}_{i\in I}$ is an open covering of $C$
such that $p$ is only included in $U_1$, then $\o_{C^{n)}}(U_i) =
\o_{C}(U_i) \text{ if } i\ne 1$ and
$$\o_{C^{n)}}(U_1) = \o_{C}(U_1) \underset {\mod
{{\o}_p}/{m_p^n}} \times k(p)\, ,$$\noindent $k(p) \longrightarrow
\mod {{\o}_p}/{m_p^n}$ being the morphism induced by the separable
elements of $\mod {{\o}_p}/{m_p^n}$ over $k(p)$, that is well-defined
because $k$ is perfect.
\end{itemize}

 $C$ and $C^{n)}$ are birationally equivalent curves, and $C^{n)}$
is a singular curve, with a unique cuspidal singularity of
multiplicity $n$ (\cite{Pa1}).

        Let ${\pi}_n
\colon C \to C^{n)}$ be the normalization morphism. If ${\bar
p}_n$ is the singular point of $C^{n)}$, one has that ${\pi}_n(p)
= {{\bar p}_n}$. We also denote by $D_n$ the fiber of ${\bar p}_n$
by ${\pi}_n$, and  ${\o}_{D_n}\simeq \mod {{\o}_p}/{m_p^n}$. Moreover,
$k({\bar p}_n) \simeq k(p)$ and ${\o}_{{\bar p}_n} \hookrightarrow
{\o}_p$.

\begin{lem} Let $B$ be a $k$-algebra. One has that
$$\mod {({\o}_{D_n} \underset {k} \otimes B)^*}/{(k(p)
\underset {k} \otimes B)^*} \simeq W^{n-1}_{k(p) \underset {k}
\otimes B}\, ,$$\noindent and, when the characteristic of $k$ is
$0$, there exists an isomorphism of groups
$$\mod {({\o}_{D_n} \underset {k}  \otimes B)^*}/{(k(p) \underset {k}
\otimes B)^*} \simeq (\mod {m_p}/{m_p^n}) \underset {k} \otimes
B$$
\end{lem}

\begin{proof} It follows from the group structure of $W^{n-1}_{k(p) \underset {k} \otimes
B}$ that the map $$\aligned W^{n-1}_{k(p) \underset {k} \otimes B}
\longrightarrow &\mod {({\o}_{D_n} \underset {k} \otimes
B)^*}/{(k(p) \underset {k} \otimes B)^*} \\
(\sum_{j_h = 1_h}^{\lambda_h} x_{j_h}\otimes b_{j_h})_h
\longmapsto &1\otimes 1 + \sum_{h=1}^{n-1} [({\bar t}_p^h\otimes
1)\cdot (\sum_{j_h = 1_h}^{\lambda_h} x_{j_h}\otimes b_{j_h})]
\endaligned$$\noindent is an isomorphism of groups, where $t_p$ is
a local parameter at $p$, ${\bar t}_p\in \mod {{\o}_p}/{m_p^n}$,
$x_{j_h}\in k(p)$ and $b_{j_h}\in B$.

        Furthermore, if the characteristic of $k$ is
$0$, the exponential map induces an isomorphism of groups between
$(\mod {m_p}/{m_p^n}) \underset {k}  \otimes B$ and $\mod
{({\o}_{D_n} \underset {k} \otimes B)^*}/{(k(p) \underset {k} \otimes
B)^*}$, which is the second part of the statement.
\end{proof}

\begin{prop} \label{p:sky}  Let ${\pi}'_{n,S} \colon D_n \times S
\to {\bar p}_n \times S$ be the morphism induced by the
desingularization map and let
 ${\phi}'_{n,2} \colon {\bar p}_n \times S
\to S$ be the natural projection. If $\overline {{{\xi}'_S}}$ is
the cokernel sheaf of the injective morphism of sheaves of abelian
groups on ${\bar p}_n \times S$ induced by ${\pi}'_{n,S}$:
$${\o}^*_{{{\bar p}_n}\times S} \to (\pi'_{n,S})_* {\o}^*_{{D_n}\times
S}\, ,$$
 and one defines ${\mathcal
W}^{n-1}_{k(p)\otimes {\o}_S}$ as the sheaf of groups on a $k$-scheme
$S$ such that it assigns $$W^{n-1}_{\Gamma (U,k(p)\otimes {\o}_S)}$$
to each open set $U$ of $S$, one has an isomorphism of sheaves of
abelian groups on $S$
$$({\phi}'_{n,2})_* \overline {{{\xi}'_S}} \simeq {\mathcal
W}^{n-1}_{k(p)\otimes {\o}_S}\, .$$
\end{prop}

\begin{proof} Since ${\phi}'_{n,2}$
is a finite morphism, one has that $R^1({\phi}'_{n,2})_*
{\o}^*_{{\bar p}_n \times S} = 0$ (\cite{Gr1}, Proposition 5.5) and,
thus, from the exact sequence of sheaves of groups on ${{\bar p}_n
\times S}$
$$0\to
{\o}^*_{{{\bar p}_n}\times S} \to (\pi'_{n,S})_* {\o}^*_{{D_n}\times S}
\to \overline {{{\xi}'_S}} \to 0$$ we have the exact sequence of
sheaves of groups on $S$
 $$0\to ({\phi}'_{n,2})_*{\o}^*_{{{\bar p}_n}\times S}
\to ({\phi}'_{n,2})_*(\pi'_{n,S})_* {\o}^*_{{D_n}\times S} \to
({\phi}'_{n,2})_* \overline {{{\xi}'_S}} \to 0.$$

    Bearing in mind that $({\phi}'_{n,2})_*{\o}^*_{{{\bar p}_n}\times
S}\simeq (k(p)\underset {k} \otimes {\o}_S)^*$ and
$$({\phi}'_{n,2})_*(\pi'_{n,S})_* {\o}^*_{{D_n}\times S}\simeq (\mod
{{\o}_p}/{m_p^n} \underset {k} \otimes {\o}_S)^*\, ,$$\noindent one has
that $({\phi}'_{n,2})_* \overline {{{\xi}'_S}}$ is the sheaf
associated with the presheaf:
$$U \longrightarrow \mod {(\mod {{\o}_p}/{m_p^n} \underset {k}
\otimes {\o}_S(U))^*}/{(k(p)\underset {k} \otimes {\o}_S(U))^*}\,
,$$\noindent where $U$ is an open set of a $k$-scheme $S$.
        Hence, it follows from the previous lemma that there exists a morphism
of sheaves between $({\phi}'_{n,2})_* \overline {{{\xi}'_S}}$ and
${\mathcal W}^{n-1}_{k(p)\otimes {\o}_S}$ that is clearly an
isomorphism.
\end{proof}

\begin{thm} The Jacobian of the singular curve
$C^{n)}$, $J_{C^{n)}}$, is an extension of the abelian variety
$J_C$ by a unipotent, connected and commutative group scheme.
Moreover, if we denote this group scheme by $K[n]$, we have that
$K[n]$ is isomorphic to the Witt scheme $W^{n-1}_{\sp S^{\cdot}
(k(p)^{\vee})}$ (Example \ref{e:witt}) and, when the
characteristic of $k$ is $0$, there exists an isomorphism of group
schemes $K[n]\simeq {\mathbb G}_a^{d(n-1)}$, where $d = {\text{
dim}}_k k(p)$.
\end{thm}

\begin{proof} Let $S$ be an arbitrary $k$-scheme. We use
$\phi_{n,2}\colon C^{n)}\times S \to S$ to denote the natural
projection and $\pi_{n,S}\colon C\times S \to C^{n)}\times S$ to
denote the morphism induced by the normalization.

        Thus, from the exact sequence of sheaves of groups on $C^{n)}\times S$
$$0\to {{\o}^*}_{C^{n)}\times S} \to (\pi_{n,S})_* {\o}^*_{C\times S}
 \to \overline {{{\xi}_S}}\to 0\, ,$$\noindent
where $\overline {{{\xi}_S}}$ is the corresponding cokernel sheaf,
and bearing in mind that $R^1 (\phi_{n,2})_* \overline
{{{\xi}_S}}=0$ because $\overline {{{\xi}_S}}$ is a concentrated
sheaf at ${\bar p}_n \times S$ which is a closed subscheme of
$C^{n)}\times S$, we obtain the morphisms of groups

\begin{equation}\label{e:chi} \xymatrix@=5pt{ 0 \ar[r] & {\Gamma} (S,(\phi_{n,2})_*
{\overline {{\xi}_S}}) \ar[r] & \Gamma (S,R^1(\phi_{n,2})_*
{\o}^*_{C^{n)}\times S}) \ar[rr]^(0.45){\pi_{n,S}^*} \ar@{=}[d]& &
\Gamma (S,R^1(\phi_{n,2})_*({\pi}_{n,S})_* {\o}^*_{C\times S})
\ar@{=}[d]
\\ & & \Pic (C^{n)})^{\bullet}(S)&  & \Pic (C)^{\bullet} (S)\quad \,
,}\end{equation} \noindent
with $\Coker \pi_{n,S}^*
\hookrightarrow H^1 (S, (\phi_{n,2})_* \overline {{{\xi}_S}})$.

    Hence, if $\overline {{{\xi}'_S}}$ is the sheaf on ${\bar p}_n \times S$ defined above,
from the commutative diagram of morphisms of schemes

$$\xymatrix@=9pt{C\times S \ar[rr]^{\pi_{n,S}} &  & C^{n)}\times S
\ar[rrd]^{\phi_{n,2}} & &
\\ & & & & S\, , \\ D_n \times S \ar[rr]^{\pi'_{n,S}} \ar@{^(->}[uu]  &  & {\bar p}_n \times S
\ar@{^(->}[uu]  \ar[rru]_{{\phi}'_{n,2}}& & }$$\noindent one has
that $(\phi_{n,2})_* \overline {{{\xi}_S}} \simeq
({\phi}'_{n,2})_* \overline {{{\xi}'_S}}$, and it follows from
Proposition \ref{p:sky} that $(\phi_{n,2})_* \overline {{{\xi}_S}}
\simeq {\mathcal W}^{n-1}_{k(p)\otimes {\o}_S}.$

        Furthermore, if $\sp B$ is an affine $k$-scheme, the sheaf ${\mathcal
W}^{n-1}_{k(p)\otimes {\o}_{\sp B}}$ is acyclic, and therefore from
(\ref{e:chi}) there exists an exact sequence of groups
$$0\to \Gamma (\sp B, {\mathcal W}^{n-1}_{k(p)\otimes {\o}_{\sp
B}})\to \Pic
 (C^{n)})^{\bullet}(\sp B)
\to \Pic (C)^{\bullet} (\sp B) \to 0$$\noindent that shows that
the induced morphism between the Picard schemes, which is denoted
by ${\bar {\pi}}_n$, is surjective.

       Accordingly, by considering the restriction of the morphism
${\bar {\pi}}_n$ to the Jacobian $J_{C^{n)}}$, we obtain the
morphisms of group schemes
$$0 \to K[n] \to J_{C^{n)}} \to J_C \to 0\, ,$$ where
${K[n]}^{\bullet} (S) =  \Gamma (S, {\mathcal
W}^{n-1}_{k(p)\otimes {\o}_{S}})$. Finally, when the characteristic
of $k$ is $0$, since a unipotent, connected, and commutative group
is isomorphic to a direct product of additive groups, one
concludes that
$$K[n]\simeq {\mathbb G}_a^{d(n-1)}.$$
 \end{proof}

\begin{cor} In arbitrary characteristic, one has that
$K[n] \simeq ({\Gamma}^{n)}_+)_{k(p)}$, with
$$[({\Gamma}^{n)}_+)_{k(p)}]^{\bullet} (S) =
\left \{1 + a_1t + \dots + a_{n-1}t^{n-1} \text { where } a_h \in
\Gamma (S,{\o}_S) \underset {k} \otimes k(p) \right \}$$ \noindent
for an arbitrary $k$-scheme $S$.
\end{cor}

\begin{cor}\label{cor:1234} One has that $K[\infty] = \underset {n} \varprojlim K[n]
\simeq ({\Gamma}_+)_{k(p)}$, with
$$[({\Gamma}_+)_{k(p)}]^{\bullet} (S) =
\left \{1 + a_1t + a_2 t^2 + \dots  \text { where } a_h \in \Gamma
(S,{\o}_S) \underset {k} \otimes k(p) \right \}$$ \noindent for an
arbitrary $k$-scheme $S$.
\end{cor}

\subsection{Group Schemes $K[{\L}^{n)}]$ over
$J_{C^{n)}}$}\label{ss-23}

    For each positive integer $n\in \mathbb N$, let us again consider
the morphism of group schemes $\bar {\pi}_n \colon J_{C^{n)}}
\longrightarrow J_C$ induced by the normalization morphism
${\pi}_n \colon C \to C^{n)}$.

    Recall from \cite{Mu} that if $X$ is a commutative group $k$-scheme,
$S$ is an arbitrary $k$-scheme and $f\colon S \to X$ is an
$S$-valued point of $X$, then $T_f$, translation by $f$, denotes
the $S$-isomorphism of schemes $$(m\circ (1_X\times f), q_2)
\colon X\times S \to X\times S\, ,$$\noindent where $m$ is the
multiplication on $X$ and $q_2\colon X\times S \to S$ is the
natural projection.

    Let $\L$ be an invertible sheaf over $J_C$. If we consider the
invertible sheaf ${\L}^{n)}$ over $J_{C^{n)}}$, ${\L}^{n)} = (\bar
{\pi}_n)^* {\L}$, from \cite{Pa2} we know of the existence
 of a group scheme $K[{\L}^{n)}]$ such that its functor of points is
$$K[{\L}^{n)}]^{\bullet} (S) = \{ f\in
J_{C^{n)}}^{\bullet} (S) \text { such that } T_f^* {\mathcal
L}^{n)}_S \simeq {\mathcal L}^{n)}_S\otimes {\bar \phi}_2^*
{\mathcal M}\}\, ,$$\noindent where $S$ is an arbitrary
$k$-scheme, ${\mathcal L}^{n)}_S$ is the pullback of the sheaf
${\mathcal L}^{n)}$ to the $k$-scheme $J_{C^{n)}}\times S$, the
map $T_f \colon J_{C^{n)}}\times S \to J_{C^{n)}}\times S$ is the
translation by the $S$-valued point $f$, ${\bar \phi}_2\colon
J_{C^{n)}}\times S \to S$ is the natural projection, and $\mathcal
M$ is an invertible sheaf over $S$.

    Furthermore, if $K[{\L}]$ is the group scheme defined by
D. Mumford in (\cite{Mu}), there exists an exact sequence of group
schemes:
$$0 \to K[n] \to K[{\L}^{n)}] \to K[{\L}] \to 0.$$

\section{Group Schemes $K[{\L}]$ over Picard Schemes of Curves}

         Let us again consider a complete, irreducible and non-singular curve $C$ over a perfect field $k$ and with a rational
point. If $p\in C$ is an arbitrary closed point, we can construct
the cuspidal curve $C^{n)}$ associated with $p$ that was
considered in the previous section.

       The aim of this section is to construct group schemes
$K[{\L}]$, associated with invertible sheaves over Picard schemes
of curves, in order to generalize the construction made by D.
Mumford for abelian varieties. Indeed, we shall provide
generalized theta groups for the schemes $\Pic (C)$ and $\Pic
(C^{n)})$, and we shall study the relation between them.

\begin{rem}  Let $X$ be a $k$-scheme such that
$X = \underset {\lambda} \coprod V_{\lambda}$, $V_{\lambda}$ being
open subschemes of $X$. With these conditions, one has an
isomorphism between the Picard groups
$$\aligned \Pic (X) &\overset {\sim} \longrightarrow
\underset {\lambda}  \prod \Pic (V_{\lambda}) \\ {\L}
 &\longmapsto ({{\L}\vert}_{V_{\lambda}})
\endaligned\, .$$
    Hence, an invertible sheaf over $X$ is determined by its restrictions
to the open subschemes $V_{\lambda}$.
\end{rem}

\subsection{Group Schemes $K[{\tilde {\L}}]$ over $\Pic C$}

         Let ${\o} (\alpha p)$ be the invertible sheaf over $C$ associated with the Weil divisor
$D = \alpha p$ for all $\alpha \in {\mathbb Z}$.

        If we consider the sheaf ${\o} (\alpha p)$ as a rational point of
the group scheme $\Pic C$, we have an automorphism of group
schemes $T_{{\o} (\alpha p)} \colon \Pic C \overset {\sim}
\longrightarrow \Pic C$, and, in particular, we have that $T_{{\o}
(\alpha p)} \colon J_{- \alpha d} \overset {\sim} \longrightarrow
J_0$ with $d = \deg(p) = {\text {dim}}_{k} k(p)$, and $J_{- \alpha
d}$ being the Jacobian that parametrizes invertible sheaves of
degree $- \alpha d$ on $C$. Let us denote by ${\mu}_{\alpha}$ the
isomorphism $T_{{\o} (-\alpha p)} \colon J_{\alpha d} \overset {\sim}
\longrightarrow J_{0}$.

        Moreover, since $J_0 = J_C \hookrightarrow \Pic (C)$
is an open subgroup (\cite{Gr2}, page 2), one has that
$\{J_{\beta}\}_{\beta \in \mathbb Z}$ is an open covering of $\Pic
(C)$. Hence, an invertible sheaf over $\Pic C$ is determined by
its restrictions to the varieties $J_{\beta}$.

\begin{defn} We shall use $Y_{\{p\}} = \underset {\alpha \in
\mathbb Z} \coprod J_{\alpha d}$ to refer to the group subscheme
of the Picard scheme $\Pic (C)$, obtained by translation of the
Jacobian $J_0$ by all the points ${\o} (\alpha p)$ with $\alpha$ an
integer.
\end{defn}

        If $p$ is a rational point, it is clear that
$Y_{\{p\}} = \Pic (C)$.

        Let $\L$ be an invertible sheaf over $J_C$. From $\L$ we
can construct the invertible sheaves ${\L}^{\alpha} =
{\mu}_{\alpha}^* {\L}$ over $J_{\alpha d}$, and we can consider
the unique invertible sheaf over $Y_{\{p\}}$, $\tilde {{\L}}$,
 whose restriction to each $J_{\alpha d}$
coincides with ${\L}^{\alpha}$. Indeed, from the canonical
immersion $i_0\colon J_{0} \hookrightarrow Y_{\{p\}}$, we have
that $i_0^* \tilde {{\L}} \simeq  {\L}$ and, if we set $i_{\alpha}
\colon J_{\alpha d} \hookrightarrow Y_{\{p\}}$ to denote the
corresponding immersion, we have that $i_{\alpha}^*  \tilde {{\L}}
\simeq {\L}^{\alpha}$.

    Let ${\mathcal C}_{k-sch.}$ be the category of algebraic $k$-schemes and
let ${\mathcal C}_{gr.}$ be the category of groups. If $S$ is an
arbitrary $k$-scheme, we can consider the functor
$$\tilde {F} \colon {\mathcal C}_{k-sch.} \longrightarrow
{\mathcal C}_{gr.}\, ,$$\noindent where $\tilde {F} (S) = \{g\in
Y_{\{p\}}^{\bullet} (S) \text { such that } T_g^* {\tilde {\L}}_S
\simeq {\tilde {\L}}_S \otimes {\bar \phi}^*_S {\M} \}$, $\M$
being an invertible sheaf over $S$, ${\bar \phi}_S \colon
Y_{\{p\}} \times S \longrightarrow S$ and ${\bar \phi}_1 \colon
Y_{\{p\}} \times S \longrightarrow Y_{\{p\}}$ being the respective
projections, and ${\tilde {\L}}_S = {\bar \phi}^*_1 {\tilde
{\L}}$.

\begin{thm}\label{t:repr} The functor $\tilde {F}$ is
representable in the category of $k$-schemes and its representant,
$K[\tilde {{\L}}]$, is a group subscheme of $\Pic (C)$.
\end{thm}

\begin{proof} If $\{J_{\beta d}\}_{\beta \in \mathbb
Z}$ is the covering of  $Y_{\{p\}}$ by the open subschemes defined
above, to prove the representability of $\tilde {F}$ is sufficient
to see that the functors $\tilde {F} \underset
{Y_{\{p\}}^{\bullet}}
 \times J_{\beta d}^{\bullet}$ are representable.

        Moreover, from the isomorphisms of varieties
${\mu}_{\beta} \colon J_{\beta d} \overset {\sim}  \longrightarrow
J_0$, one has that $$\tilde {F} \underset {Y_{\{p\}}^{\bullet}}
 \times J_{\beta d}^{\bullet} \simeq
\left ( \tilde  {F} \underset {Y_{\{p\}}^{\bullet}} \times
J_0^{\bullet} \right ) \underset {J_0^{\bullet}} \times J_{\beta
d}^{\bullet}\ .$$

Hence, to prove the representability of $\tilde {F}$, one only has
to see that $\tilde {F} \underset {Y_{\{p\}}^{\bullet}}
 \times J_0^{\bullet}$ is representable.

        Furthermore, since ${\L}$
is an invertible sheaf over $J_0$, which is an abelian variety,
there exists a group scheme $K[{\L}] \subseteq J_0$ \cite{Mu},
 where $$K[{\L}]^{\bullet} (S) = \{f\in J_0^{\bullet} (S)
\text { such that } T_f^* {\L}_S \simeq {\L}_S \otimes {\phi_2}^*
{\M} \}$$ for a certain invertible sheaf $\M$ over each $k$-scheme
$S$, $\phi_2 \colon J_0 \times S \to S$ being the natural
projection.

    Thus, to see the representability of  $\tilde {F}$, we only need
to prove that $$\tilde {F} \underset {Y_{\{p\}}^{\bullet}} \times
J_0^{\bullet} \simeq K[{\L}]^{\bullet}\, .$$

    Accordingly, given an arbitrary $k$-scheme $S$, and setting
$$(i_0)_S \colon J_0 \times S \hookrightarrow Y_{\{p\}}\times S$$\noindent
to denote the morphism of schemes induced by $i_0$, one has that
$(i_0)_S^* {\tilde {\L}}_S \simeq {\L}_S$. If we also consider
$f\in (\tilde {F} \underset {Y_{\{p\}}^{\bullet}} \times
J_0^{\bullet})(S)$, we deduce that $(i_0)_S \circ T_f = T_{\bar f}
\circ (i_0)_S$, with ${\bar f} \in {\tilde {F}}(S)$, and it
follows from this that there exists a natural injective morphism
$\tilde {F} \underset {Y_{\{p\}}^{\bullet}} \times J_0^{\bullet}
\overset {j} \hookrightarrow K[{\L}]^{\bullet}$.

        To conclude, we have to see that the map $j$ is surjective.
Let us now consider $g\in  K[{\L}]^{\bullet} (S)$ for a certain
$k$-scheme $S$. One has that  $$T_g^* {\L}_S \simeq {\L}_S \otimes
\phi_2^* {\mathcal M}\, .$$
        If $({\mu}_{\beta})_S
\colon J_{\beta d} \times S \overset {\sim} \longrightarrow J_0
\times S$ are the morphisms induced by ${\mu}_{\beta}$, it follows
from the isomorphism of sheaves
 $$T_g^* ({\mu}_{\beta})_S^* {\L}_S \simeq
 ({\mu}_{\beta})_S^* T_g^* {\L}_S$$
that $T_g^* {\L}^{\beta}_S \simeq {\L}^{\beta}_S \otimes
(\phi_2^{\beta})^* {\mathcal M}$ with
 ${\L}^{\beta}_S = ({\mu}_{\beta})_S^* {\L}_S$,
 $\phi_2^{\beta} \colon J_{\beta d} \times S \to S$
 being the natural projection.

        Finally, since ${\tilde {\L}}_S = \underset {\beta \in
\mathbb Z} \coprod {\L}^{\beta}_S$, from the commutative diagrams
of morphisms of schemes
$$\xymatrix@=8pt{Y_{\{p\}}\times S
\ar[rrd]^{{\bar \phi}_S} & &
\\ & & S\, , \\ J_{\beta d} \times S
\ar@{_(->}[uu]  \ar[rru]_{\phi_2^{\beta}}& & }$$\noindent one has
that $T_g^* {\tilde {\L}}_S \simeq {\tilde {\L}}_S \otimes {\bar
\phi}^*_S {\M}$, and, therefore, $g\in (\tilde {F} \underset
{Y_{\{p\}}^{\bullet}} \times J_0^{\bullet}) (S)$.

        Hence $\tilde {F}$ is representable in the category of
$k$-schemes. Moreover, its representant $K[\tilde {{\L}}]$ is a
group subscheme of $Y_{\{p\}}$ and is hence a group subscheme of
$\Pic (C)$.  \end{proof}

\begin{cor} If $\L$ is an invertible sheaf over
$J_C$ such that $K[{\L}]$ is reduced, then the group scheme
$K[\tilde  {{\L}}]$ is reduced.
\end{cor}

\begin{proof} The statement is a direct consequence of the fact
that the connected component of the identity element of $K[\tilde
{{\L}}]$ is $K[{{\L}}]$ which, with the hypothesis of the
corollary, is a reduced $k$-scheme.
\end{proof}

\subsection{Group Schemes $K[{\tilde {{\L}}}^{n)}]$ over $\Pic (C^{n)})$}

   In the previous section we set
$Y_{\{p\}} = \underset {\alpha \in \mathbb Z} \coprod J_{\alpha
d}$ to denote the group subscheme of $\Pic (C)$ obtained by
translation of the Jacobian $J_0$ by all the points ${\o} (\alpha
p)$.

        For each positive integer $n\in \mathbb N$, let us again consider the morphism
of group schemes $\bar {\pi}_n \colon \Pic (C^{n)})
\longrightarrow \Pic C$ induced by the normalization morphism
${\pi}_n \colon C \to C^{n)}$.

\begin{defn} We shall use $Y_{\{p\}}^{n}$ to denote the group subscheme of
$\Pic (C^{n)})$ which is the fiber of $Y_{\{p\}}$ by $\bar
 {\pi}_n$;
 i.e. $Y_{\{p\}}^{n} = {\bar {\pi}_n}^{-1} (Y_{\{p\}})$.
\end{defn}

        If we set $(Y^n_{\{p\}})_{\beta} = {\bar
{\pi}_n}^{-1} (J_{\beta})$, and $d = \deg(p) = {\text {dim}}_{k}
k(p)$, we have that  $$Y_{\{p\}}^{n} = \underset {\alpha \in
\mathbb Z} \coprod (Y^n_{\{p\}})_{\alpha d}\, ,$$\noindent and,
similar to above, when $p$ is a rational point, one has that
$Y_{\{p\}}^{n} = \Pic (C^{n)})$.

        In particular, we know that $(Y_{\{p\}}^n)_0 = {\bar {\pi}_n}^{-1} (J_C) =
J_{C^{n)}}$ is a group scheme that is an extension of $J_C$ by a
unipotent, connected and commutative group scheme $K[n] = \Ker
\bar {\pi}_n$.

        Moreover, it follows from translating by a fixed rational point of $\Pic
(C^{n)})$ in the fiber of each invertible sheaf ${\o} (\alpha p)$,
that there exist isomorphisms of schemes
$${\mu}_{\beta}^n \colon (Y_{p}^n)_{\beta d} \overset {\sim}
\longrightarrow J_{C^{n)}}\, .$$

    Let us again consider an invertible sheaf $\L$ over $J_C$, and
let us denote by $\tilde {\L}$ its extension to $Y_{\{p\}}$,
defined previously. If we set ${\tilde {{\L}}}^{n)} = (\bar
{\pi}_n)^* {\tilde {\L}}$ and ${\L}^{n)} = (\bar {\pi}_n)^* {\L}$,
keeping the above notations, by construction we have that
${\L}^{n)}$ is an invertible sheaf over $J_{C^{n)}}$. Furthermore,
if $(i_{\alpha}^{n})\colon (Y^n_{\{p\}})_{\alpha d}
\hookrightarrow \Pic (C^{n)})$ is the natural immersion, one has
that $(i_{\alpha}^{n})^* {\tilde {{\L}}}^{n)} \simeq
({\L}^{\alpha})^{n)}$ for each $\alpha \in \mathbb Z$, where
$$({\L}^{\alpha})^{n)} =
 (\bar {\pi}_n)^* {\L}^{\alpha}
= ({\mu}^n_{\alpha})^*{\L}^{n)}\, .$$

    Recall from Subsection \ref{ss-23} that, according to the
results of \cite{Pa2}, there exists a group scheme
$K[{\L}^{n)}]\subseteq J_{C^{n)}}$, whose functor of points is
$$K[{\L}^{n)}]^{\bullet} (S) = \{ f\in
J_{C^{n)}}^{\bullet} (S) \text { such that } T_f^* {\mathcal
L}^{n)}_S \simeq {\mathcal L}^{n)}_S\otimes {\bar \phi}_2^*
{\mathcal M}\}\, .$$

    Let $S$ be an arbitrary $k$-scheme.
Analogously to the previous subsection, we define the functor
$$\tilde {F_n} \colon {\mathcal C}_{k-sch.}
\longrightarrow {\mathcal C}_{gr.}\, , $$\noindent where $\tilde
{F_n} (S) = \{f\in (Y_{\{p\}}^n)^{\bullet} (S) \text { such that }
T_f^* {\tilde {{\L}}}^{n)}_S \simeq {\tilde {{\L}}}^{n)}_S \otimes
{\bar \phi}^*_{n,S} {\mathcal M} \}$, $\mathcal M$ is an
invertible sheaf over $S$,
 ${\bar \phi}_{n,S} \colon Y_{\{p\}}^n \times S
\longrightarrow S$ and ${\bar \phi}_{n,1} \colon Y_{\{p\}}^n
\times S \longrightarrow Y_{\{p\}}^n$ are the respective
projections, and ${\tilde {{\L}}}^{n)}_S = {\bar \phi}^*_{n,1}
{\tilde {{\L}}}^{n)}$.

\begin{thm} For each integer $n>1$, the functors $\tilde {F_n}$ are representable
in the category of $k$-schemes.

     Their respective representants are group schemes $K[{\tilde {{\L}}}^{n)}] \hookrightarrow \Pic (C^{n)})$.
\end{thm}

\begin{proof}  Fixing a positive integer $n>1$, let us consider the covering of $Y_{\{p\}}^n$
by the open subschemes $\{(Y^n_{\{p\}})_{\beta d}\}_{\beta \in
\mathbb Z}$, such that $(Y^n_{\{p\}})_{\beta d}$ is the fiber of
the Jacobian $J_{\beta d}$ by the morphism $\bar {\pi}_n$.

From the isomorphisms ${\mu}^n_{\beta} \colon (Y^n_{\{p\}})_{\beta
d} \overset {\sim} \longrightarrow J_{C^{n)}}$, and with similar
arguments to the proof of Theorem \ref{t:repr}, one has that
$${\tilde {F_n}} \underset
{(Y_{\{p\}}^n)^{\bullet}} \times J_{C^{n)}}^{\bullet} \simeq
K[{\L}^{n)}]^{\bullet}$$ and, hence,  $${\tilde {F_n}} \underset
{(Y_{\{p\}}^n)^{\bullet}} \times (Y^n_{\{p\}})_{\beta d}^{\bullet}
\simeq K[{\L}^{n)}]^{\bullet} \underset {J_{C^{n)}}^{\bullet}}
\times (Y_{\{p\}}^n)_{\beta d}^{\bullet}$$
    Therefore $\tilde {F_n}$ is representable because the functors ${\tilde {F_n}}
\underset {(Y_{\{p\}}^n)^{\bullet}} \times (Y^n_{\{p\}})_{\beta
d}^{\bullet}$ are representable and $\{(Y^n_{\{p\}})_{\beta
d}\}_{\beta \in \mathbb Z}$ is an open covering of $Y_{\{p\}}^n$
(\cite{EGA}). We shall use
 $K[{\tilde {{\L}}}^{n)}]$ to denote its representant, which is a group
subscheme of $\Pic (C^{n)})$.
\end{proof}

\begin{cor} If $\L$ is an invertible sheaf over
$J_C$ such that $K[\L]$ is reduced, one has that the group scheme
 $K[{\tilde {{\L}}}^{n)}]$ is reduced for all $n>1$.
\end{cor}

\begin{proof} It is known that the group scheme
 $K[{\L}^{n)}]$ is reduced (\cite{Pa2}) and therefore we are done,
 bearing in mind that this group scheme is the connected
 component of the identity element of the $k$-scheme $K[{\tilde {{\L}}}^{n)}]$.
\end{proof}

        Our purpose now is to determine the relation between the group
scheme $K[{\tilde {\L}}]$, defined in the above subsection, and
the group schemes $K[{\tilde {{\L}}}^{n)}]$.

\begin{lem}\label{l:sep} For each separated $k$-scheme
$S$, if ${\bar {\pi}}_{n,S}\colon Y_{\{p\}}^{n)} \times S \to
Y_{\{p\}} \times S$ is the morphism of group schemes induced by
${\bar {\pi}}_n$, then the natural morphism of groups
$${\bar {\pi}}_{n,S}^* \colon \Pic (Y_{\{p\}} \times S) \longrightarrow \Pic (Y_{\{p\}}^{n)}
\times S)$$ is injective.
\end{lem}

\begin{proof} With the above notations, one has that $$Y_{\{p\}} \times S = \underset {\alpha
\in {\mathbb Z}} \coprod \Big ( (Y_{\{p\}})_{\alpha
 d}\times S \Big )\quad \text { and }\quad Y_{\{p\}}^{n} \times S =
\underset  {\beta
\in {\mathbb Z}} \coprod \Big ( (Y_{\{p\}}^n)_{\beta
 d}\times S \Big )\ ,$$ and there exist isomorphisms of schemes $${\mu}_{\alpha} \colon (Y_{\{p\}})_{\alpha
 d} \overset {\sim} \to J_C \text{ and } {\mu}_{\beta}^{n)} \colon
(Y_{\{p\}}^{n})_{\beta  d} \overset {\sim} \to J_{C^{n)}}\, .$$

    Moreover, from \cite{Pa2} one has that the natural morphism of groups  $$\Pic (J_C \times S) \to \Pic
(J_{C^{n)}} \times S)$$\noindent is injective for each separated
$k$-scheme $S$.

    Thus, the morphism of groups ${\bar {\pi}}_{n,S}^* \colon \Pic (Y_{\{p\}} \times S) \longrightarrow
\Pic (Y_{\{p\}}^{n} \times S)$ is injective because it is in each
component.
\end{proof}

\begin{prop} If $S$ is a separated $k$-scheme, one has that
$$[K[{\tilde {\L}}]^{\bullet} \underset
{Y_{\{p\}}^{\bullet}} \times (Y_{\{p\}}^n)^{\bullet}] (S) \simeq
K[{\tilde {{\L}}}^{n)}]^{\bullet} (S)\, .$$
\end{prop}

\begin{proof}  Given an arbitrary $k$-scheme $S$,
one has that
$$[K[{\tilde {\L}}]^{\bullet} \underset
{Y_{\{p\}}^{\bullet}} \times (Y_{\{p\}}^n)^{\bullet}] (S)
 \overset {i}  \hookrightarrow K[{\tilde {{\L}}}^{n)}]^{\bullet} (S)$$
because if $f\in [K[{\tilde {\L}}]^{\bullet} \underset
{Y_{\{p\}}^{\bullet}} \times (Y_{\{p\}}^n)^{\bullet}] (S)
\subseteq (Y_{\{p\}}^n)^{\bullet} (S)$, from the commutative
diagram of morphisms of schemes
$$\xymatrix@=9pt{Y_{\{p\}}^{n)}
\times S  \ar[dd]_{{\bar {\pi}}_{n,S}} &  & & Y_{\{p\}}^{n)}
\times S \ar[lll]_{T_f} \ar[dd]^{{\bar {\pi}}_{n,S}}
\ar[rrrd]^{{\bar \phi}_{n,S}} & & &
\\ & & & & & & S\, , \\ Y_{\{p\}} \times S    & &
  & Y_{\{p\}} \times S \ar[lll]_{T_{{\bar {\pi}}_n
(f)}}
  \ar[rrru]_{{\bar \phi}_{S}}& & &}$$\noindent
there exists the sequence of isomorphisms
$$T_f^* {\tilde {{\L}}}^{n)}_S  \simeq
{\bar {\pi}}_{n,S}^* (T_{{\bar {\pi}}_n (f)}^* {\tilde {\mathcal
L}}_S) \simeq {\bar {\pi}}_{n,S}^* ({\tilde {\mathcal L}}_S
\otimes {\bar \phi}_S^* {\mathcal N})\simeq {\tilde {{\L}}}^{n)}_S
\otimes {\bar \phi}^*_{n,S} {\mathcal N}$$ for a certain
invertible sheaf ${\mathcal N}$ over $S$, and one deduces that
$f\in K[{\tilde {{\L}}}^{n)}]^{\bullet} (S)$.

     Thus, to conclude we must see that the map $i$ is surjective when $S$ is
a separated $k$-scheme.

        Let us consider $g\in K[{\tilde {{\L}}}^{n)}]^{\bullet} (S)$, such that $T_g^* {\tilde {{\L}}}^{n)}_S \simeq
{\tilde {{\L}}}^{n)}_S \otimes {\bar \phi}^*_{n,S} {\mathcal M}$,
$\mathcal M$ being an invertible sheaf over $S$.

    Bearing in mind that $T_g^* ({\bar {\pi}}_{n,S}^*
{\tilde {\mathcal L}}_S) \simeq {\bar {\pi}}_{n,S}^* (T_{{\bar
{\pi}}_n(g)}^* {\tilde {\mathcal L}}_S)$, from Lemma \ref{l:sep}
one has that $$T_{{\bar {\pi}}_n(g)}^* ({\tilde {\mathcal L}}_S)
\simeq ({\tilde {\mathcal L}}_S\otimes {\bar \phi}_S^* {\mathcal
M})\otimes {\bar \phi}_S^* {\mathcal M'} \simeq {\tilde {\mathcal
L}}_S\otimes {\bar \phi}_S^* \Bar {\mathcal M}\, ,$$\noindent with
$\Bar {\mathcal M}\simeq {\mathcal M}\otimes {\mathcal M'}$, and
${\mathcal M'}$ being an invertible sheaf over $S$.

        Then, ${\bar {\pi}}_n (g) \in K[{\tilde
{{\L}}}]^{\bullet} (S)$, i.e.  $g\in [K[{\tilde {\L}}]^{\bullet}
\underset {Y_{\{p\}}^{\bullet}} \times (Y_{\{p\}}^n)^{\bullet}]
(S)$.
 \end{proof}

\begin{prop}\label{p:ex} For each $n>1$, if $K[n] =
\Ker {\bar {\pi}_n}$, there exists an exact sequence of group
schemes
$$0 \to K[n] \to K[{\tilde {{\L}}}^{n)}]
 \overset {\bar {\pi}_n}
\longrightarrow K[{\tilde {\L}}] \to 0\, .$$

    Hence the natural morphism $K[{\tilde {{\L}}}^{n)}]
 \overset {\bar {\pi}_n}
\longrightarrow K[{\tilde {\L}}]$ is a geometric quotient of group
schemes by the action of $K[n]$, which is a unipotent, connected
and commutative group scheme.
\end{prop}

\begin{proof} From the morphisms of schemes $$0 \to K[n] \to Y_{\{p\}}^n \overset {\bar {\pi}_n}
\longrightarrow Y_{\{p\}} \to 0\, ,$$\noindent and bearing in mind
that functors of points are sheaves for the faithfully flat and
quasi-compact topology, the claim is deduced immediately from the
previous proposition.
\end{proof}

\subsection{Group Schemes $K[(\tilde
{{\L}_{\Theta}})^{n)}]$.}

    If $g$ is the genus of $C$, let us now consider the invertible sheaf ${\L}_{\Theta}^{g-1}$ associated with the
principal polarization of the Jacobian  $J_{g-1} \subseteq \Pic
C$, which is defined canonically.  Fixing a rational point in the
curve $C$, we can translate ${\L}_{\Theta}^{g-1}$ to the Jacobian
$J_C$, and we shall use ${\L}_{\Theta}$ to denote this principal
polarization.

    Moreover, if $d = \deg (p)$, and  ${\mu}_{\alpha} \colon J_{\alpha d}
\overset {\sim} \longrightarrow J_C$ is again the isomorphism of
schemes obtained by translation by ${\o} (-\alpha p)$, we shall use
$\tilde {{\L}_{\Theta}}$ to denote the unique invertible sheaf
over $Y_{\{p\}}$ which coincides with ${\L}_{\Theta}^{\alpha} =
 {\mu}_{\alpha}^* {\L}_{\Theta}$ when we restrict to
$J_{\alpha d}$.

    If $(\tilde
{{\L}_{\Theta}})^{n)} = {\bar {\pi}_n}^* \tilde {{\L}_{\Theta}}$,
bearing in mind the previous subsection, there exist group schemes
 $K[(\tilde {{\L}_{\Theta}})^{n)}]
\subseteq Y_{\{p\}}^n$, whose structure will be determined in this
subsection.

\begin{lem}\label{l:zeta} One has that: $$K[\tilde
{{\L}_{\Theta}}]^{\bullet} (\sp k) \simeq \{{\o} (\alpha p)\}_{\alpha
\in \mathbb Z} \subseteq \Pic (C)^{\bullet} (\sp k)\, ,$$\noindent
i.e., $K[\tilde {{\L}_{\Theta}}]^{\bullet} (\sp k) \simeq \mathbb
Z$.
\end{lem}

\begin{proof} From the definition of $K[\tilde
{{\L}_{\Theta}}]$, we have that $$\aligned \quad  &K[\tilde
{{\L}_{\Theta}}]^{\bullet} (\sp k) = \{ x\in Y_{\{p\}}^{\bullet}
(\sp k) \text { such that } T_x^* (\tilde {{\L}_{\Theta}}) \simeq
(\tilde {{\L}_{\Theta}}) \} \\ &= \{  x\in Y_{\{p\}}^{\bullet}
(\sp k) \text { such that } T_x^* ({\L}^{\alpha)}_{\Theta})
 \simeq {\L}_{\Theta}
\text { for a certain } \alpha \in \mathbb Z \} \\ &= \{ x\in
Y_{\{p\}}^{\bullet} (\sp k) \text { such that } T_x^*
(T^*_{{\o}(-\alpha p)} {\L}_{\Theta}) \simeq {\L}_{\Theta} \text {
for a certain } \alpha \in \mathbb Z \}\, ,\endaligned$$\noindent
where ${\L}_{\Theta}$ is the sheaf associated with the principal
polarization of $J_C$.

    Thus, since $K[{{\L}_{\Theta}}]^{\bullet}
(\sp k) = \{e\}$, one concludes that
$$K[\tilde{{\L}_{\Theta}}]^{\bullet} (\sp k) \simeq \{{\o} (\alpha
p)\}_{\alpha \in \mathbb Z}\, .$$
\end{proof}

\begin{prop} If we consider the exact sequence of groups $$0 \to K[n]^{\bullet} (\sp k) \to
K[(\tilde {{\L}_{\Theta}})^{n)}]^{\bullet} (\sp k) \overset {\bar
{\pi}_n} \longrightarrow K[\tilde {{\L}_{\Theta}}]^{\bullet} (\sp
k) \to 0\, ,$$\noindent there exists a natural section $\sigma
\colon K[\tilde {{\L}_{\Theta}}]^{\bullet} (\sp k) \to K[(\tilde
{{\L}_{\Theta}})^{n)}]^{\bullet} (\sp k)$ of the morphism ${\bar
{\pi}_n}$.\
\end{prop}

\begin{proof} If
${\pi}_n \colon C \longrightarrow C^{n)}$ is the normalization
morphism, one has that  ${\bar {\pi}_n}$ coincides with
${\pi}_n^*$ over the invertible sheaves of the singular curve
$C^{n)}$.

 It follows from Lemma \ref{l:zeta} that
$K[\tilde {{\L}_{\Theta}}]^{\bullet} (\sp k) = \{{\o}(\alpha
p)\}_{\alpha \in \mathbb Z}$, where ${\o}(\alpha p)$ is the
invertible sheaf associated with the Cartier divisor $$D_{\alpha}
= \{(t_p^{-\alpha},U_1),(1,U_i)\}_{i\ne 1}\, ,$$\noindent  $t_p$
being the generator of the maximal ideal  ${\mathfrak m}_p
\subseteq {\o}_p$, and $\{U_i\}_{i\in I}$ being a covering of $C$
such that $p$ is only included in $U_1$; hence, we can identify
${\o}(\alpha p)$ with $t_p^{-\alpha}$.

    Moreover, since ${\Sigma}^*_C = {\Sigma}^*_{C^{n)}}$, the Cartier divisor
$D_{\alpha} = \{(t_p^{-\alpha},U_1),(1,U_i)\}$ determines over
$C^{n)}$ an invertible sheaf ${\L}_{\alpha p}^{n)}$.

    Furthermore, ${\L}_{\alpha p}^{n)} (U_1) \simeq t_p^{\alpha}\cdot {\o}_{C^{n)}}
(U_1)$, ${\o}(\alpha p) (U_1) \simeq t_p^{\alpha}\cdot {\o}_{C} (U_1)$,
and for $i\ne 1$, there also exist isomorphisms ${\L}_{\alpha
p}^{n)} (U_i) \simeq {\o}_{C^{n)}} (U_i)$ and ${\o}(\alpha p) (U_i)
\simeq {\o}_{C} (U_i)$.

    By construction, it is clear that
${\o} (\alpha p)\simeq {\pi}^*_n {\L}_{\alpha p}^{n)}$, and hence
there exists a section of the morphism ${\bar {\pi}_n}$.
\end{proof}

    By considering the category of locally noetherian schemes, we
shall now prove that $K[(\tilde {{\L}_{\Theta}})^{n)}]$ is a
direct product, where one of the factors is the connected
component of its identity element.

    Let us set ${\mathbb Z}_{*} = \underset {{\alpha}\in {\mathbb Z}}
 \coprod \sp k$. For each locally noetherian scheme $S$, one has that  $${\mathbb
Z}_{*}^{\bullet} (S) = {\Apl}_{cont.} (S,{\mathbb Z})\,
,$$\noindent where
 $\mathbb Z$ is a discrete topological space. Moreover, we have
 a structure of group scheme in ${\mathbb Z}_{*}$ that coincides,
 when we consider points on connected $k$-scheme,
 with the group structure of $\mathbb Z$.

\begin{thm}\label{t:323} For each $n\in \mathbb N$, one has an isomorphism of group schemes
 $$K[(\tilde
{{\L}_{\Theta}})^{n)}] \simeq K[(\tilde {{\L}_{\Theta}})^{n)}]_{0}
 \times {\mathbb Z}_{*}\, ,$$\noindent where
$K[(\tilde {{\L}_{\Theta}})^{n)}]_{0}$ is the connected component
of the identity element of the scheme $K[(\tilde
{{\L}_{\Theta}})^{n)}]$.
\end{thm}

\begin{proof} Let us use
${\lambda}_{\alpha}^{n)}\colon K[(\tilde {{\L}_{\Theta}})^{n)}]
 \overset {\sim} \longrightarrow
K[(\tilde {{\L}_{\Theta}})^{n)}]$ to denote the isomorphism of
schemes induced by translation by the rational point
${\L}_{{\alpha}p}^{n)}$.

        Moreover, since the schemes of the statement are locally noetherian, the sets
of their connected components are locally finite, and they are
therefore locally connected (\cite{EGA}, page 50). Hence, their
connected components are open subschemes, and to prove the theorem
it is sufficient to see that there exist isomorphisms when we
consider $S$-valued points, $S$ being a connected $k$-scheme.

        For a connected $k$-scheme $S$, we consider the map  $$\varphi
\colon (K[(\tilde {{\L}_{\Theta}})^{n)}]_{0}
 \times {\mathbb Z}_{*})^{\bullet} (S)
\longrightarrow K[(\tilde {{\L}_{\Theta}})^{n)}]^{\bullet} (S)\, ,
$$\noindent defined by $\varphi (f,\alpha) = {\lambda}^{n)}_{\alpha} \circ
f$.

        Since each $S$-valued point of $K[(\tilde
{{\L}_{\Theta}})^{n)}]$, $S$ being a connected $k$-scheme, takes
values in a unique connected component, which is determined by an
integer number $\alpha$, one has that the map $\varphi$ is
bijective.

        Furthermore, it follows from the isomorphism of invertible
sheaves ${\L}_{{\alpha}p}^{n)}\otimes
 {\L}_{{\beta}p}^{n)}
\simeq {\L}_{({\alpha + \beta})p}^{n)}$,  that $\varphi$ is a
morphism of groups for each locally noetherian $k$-scheme $S$, and
one concludes that $K[(\tilde {{\L}_{\Theta}})^{n)}]$ and
$K[(\tilde {{\L}_{\Theta}})^{n)}]_{0} \times {\mathbb Z}_{*}$ are
isomorphic as group schemes.
\end{proof}

\begin{cor}\label{c:char} One has that  $$K[(\tilde
{{\L}_{\Theta}})^{n)}] \simeq K[n]
 \times {\mathbb Z}_{*}.$$
\end{cor}

\begin{proof}
    Bearing in mind that $K[(\tilde {{\L}_{\Theta}})]_0 =
K[({\L}_{\Theta})] = \{e\}$, from Proposition \ref{p:ex} we have
the existence of an isomorphism of groups $K[(\tilde
{{\L}_{\Theta}})^{n)}]_0 \simeq K[n]$, and the claim is deduced.
\end{proof}

\section{Generalized Local Symbols}

\subsection{The Contou-Carr\`{e}re symbol associated with a rational
point}\label{ss:41}

Let $C$ be a complete, irreducible and non-singular curve over a
perfect field $k$ and let $p\in C$ be a rational point. We set:
$$K^{\Theta}_{\{p\}} = \underset {n} \varprojlim
K[(\tilde {{\mathcal L}_{\Theta}})^{n)}]\, ,$$\noindent $K[(\tilde
{{\mathcal L}_{\Theta}})^{n)}]$ being the group scheme constructed
in the previous section.

        Bearing in mind that $K[(\tilde
{{\mathcal L}_{\Theta}})^{n)}] \simeq K[n] \times {\Bbb Z}_{*}$
(Corollary \ref{c:char}),  where $K[n]$ is the kernel of the
natural morphism of group schemes $\Pic (C^{n)}) \to \Pic (C)$,
from Corollary 2.5 one has that
$$K^{\Theta}_{\{p\}} \simeq {\Bbb Z}_{*} \times K[\infty] \simeq {\Bbb Z}_{*}
\times {\Gamma}_+\, ,$$ \noindent  ${\Gamma}_+$ being the
representant of the functor on groups:
$$S\rightsquigarrow {\Gamma}_{+}(S) = \left \{ {\begin{matrix} \text
{series } 1 + a_1z + a_2z^2 + \dots \\ \text {where } a_i\in
H^0(S, {\o}_S) \end{matrix}} \right \}.$$

        Moreover, the groups schemes ${\Bbb Z}_{*}$ and ${\Bbb G}_m$ are
autodual and the universal character coincides, for each connected
$k$-scheme $S$, with the morphism of groups
$$\aligned {\Bbb Z}_{*}^{\bullet} (S) \times {\Bbb G}_m^{\bullet}
(S) &\longrightarrow {\Bbb G}_m^{\bullet} (S) \\ (\alpha ,\lambda)
&\longmapsto {\lambda}^{\alpha}\quad \, . \endaligned$$

    Let us also consider the formal group scheme ${\Gamma}_{-}$,
whose functor of points is
$$S\rightsquigarrow {\Gamma}_{-}(S) = \left \{
{\begin{matrix} \text {series } a_nz^{-n}+ \dots + a_1z^{-1} + 1
\\ \text {where } a_i\in H^0(S, {\o}_S) \text { are } \\ \text {
nilpotents and n is arbitrary }\end{matrix}} \right \}\, .$$

        It follows from the duality between the group
schemes ${\Gamma}_+$ and ${\Gamma}_{-}$ (\cite{Haz}, page 500)
that ${\hat K}^{\Theta}_{\{p\}} \simeq {\Bbb G}_m \times
{\Gamma}_{-}$, ${\hat K}^{\Theta}_{\{p\}}$ being the Cartier dual
of the group scheme $K^{\Theta}_{\{p\}}$.

        Let us now set ${\Gamma}_p =
K^{\Theta}_{\{p\}} \times {\hat K}^{\Theta}_{\{p\}}$. We have that
$${\Gamma}_p \simeq {\Bbb Z}_{*}\times {\Gamma}_{+}\times
 {\Bbb G}_m\times
{\Gamma}_{-}\, ,$$\noindent and therefore ${\Gamma}_p$ is a
locally connected group scheme.

        If $\hat {{\o}_p}$ is the completion of the local ring ${{\o}_p}$ and
$(\hat {{\o}_p})_{(0)}$ is the field of fractions of $\hat {{\o}_p}$,
for each connected $k$-scheme $S$ one has that
$${\Gamma}_p^{\bullet} (S) \simeq  H^{0}(S,{\o}_S) ((z))^{*}\, ,$$
and hence ${\Gamma}_p^{\bullet} (\sp k) \simeq k((z))^{*} \simeq
(\hat {{\o}_p})^{*}_{(0)}$, which implies that
$${\Sigma}^{*}_C \hookrightarrow {\Gamma}_p^{\bullet} (\sp
k)\, .$$

\begin{defn} We shall use the term Heisenberg group scheme
associated with ${\Gamma}_p$, denoting this by ${\mathcal H}
({\Gamma}_p)$, to refer to the scheme $${\mathcal H} ({\Gamma}_p)
= {\Bbb G}_m \times K^{\Theta}_{\{p\}} \times {\hat
K}^{\Theta}_{\{p\}}$$ together with the group law $$(\alpha, x,
l)\cdot ({\alpha}', x',l') = (\alpha \cdot {\alpha}' \cdot l(x'),
x\cdot x',l\cdot l')$$ \noindent for $S$-valued points, $S$ being
a $k$-scheme.
\end{defn}

        If $e_{{\mathcal H} ({\Gamma}_p)}
(x,y)$ is the commutator in the Heisenberg group scheme associated
with ${\Gamma}_p$, and since ${\Gamma}_p\simeq {\hat {\Gamma}}_p$,
the map
$$\aligned \varphi \colon {\Gamma}_p &\longrightarrow {\hat
{\Gamma}}_p \\ x &\longmapsto e_{{\mathcal H} ({\Gamma}_p)}
(x,\quad )
\endaligned$$ is an isomorphism of groups, because if $x = (a,b)\in
K^{\Theta}_{\{p\}} \times {\hat K}^{\Theta}_{\{p\}}$, then
$$\varphi (a,b) = (a^{-1},b)\, .$$

        Therefore, ${\mathcal H} ({\Gamma}_p)$ satisfies the characterization
of a Heisenberg group of an extension by the multiplicative group
(\cite{MNN}, page 2).

    Furthermore, since $e_{{\mathcal H}
({\Gamma}_p)}$ is a 2-cocycle, it determines an element of the
cohomology group $H_{reg}^2 ({\Gamma}_p,{\mathbb G}_m)$, which
contains the classes of 2-cocycles that are morphisms of schemes.
We shall now recall from \cite{Pa3} the definition, from this
cohomology class, of the Contou-Carr\`{e}re symbol as a morphism
of schemes.

\begin{defn} \label{d:valuation} If $S$ is a connected $k$-scheme and $f\in
{\Gamma}_p^{\bullet} (S) \simeq {\mathbb Z}_{*}^{\bullet} (S)
\times ({\Gamma}_{+}\times
 {\Bbb G}_m\times
{\Gamma}_{-})^{\bullet} (S)$, we shall call its component in
${\mathbb Z}_{*}^{\bullet} (S)$, which is an integer number,
$v(f)$, ``the valuation of $f$''.
\end{defn}
    Bearing in mind that ${\Gamma}_p = {\tilde {\Gamma}}$,
 ${\tilde {\Gamma}}$ being the formal group scheme studied in
\cite{Pa3}, it follows from \cite{Pa3}, Theorem 3.5, that:

\begin{thm} There exists a unique element $(\quad,
\quad)_{p}$ in the cohomology class $[e_{{\mathcal H}
 ({\Gamma}_p)}] \in H_{reg}^2
({\Gamma}_p,{\mathbb G}_m)$ satisfying the conditions:

\begin{itemize}
\item $(f,g\cdot g')_p = (f,g)_{p} \cdot (f,g')_{p}$;
\item $(f,g)_{p} = e_{{\Bar {\mathcal H}}
({\tilde {{\Gamma}_{p}}})} (f,g)$ if  $v(f)= 0$;
\item $(f,-f)_{p} = 1$
\end{itemize}
\noindent for $f,g,g' \in {{\Gamma}_p}^{\bullet} (S)$, with $S$ a
connected $k$-scheme.
\end{thm}

\begin{cor}\label{c:2003} If $S$ is a connected
$k$-scheme and $u, w \in {{\Gamma}}_p^{\bullet}  (S)$ with
$$\begin{aligned} \quad &u = \lambda z^n{\prod}_{i=1}^{l}
(1-a_{-i}z^{-i}) {\prod}_{i=1}^{\infty} (1-a_{i}z^{i}) \\
&w = \mu z^m {\prod}_{j=1}^{h} (1-b_{-j}z^{-j})
{\prod}_{j=1}^{\infty} (1-b_{j}z^{j})\, ,
\end{aligned}$$\noindent where $\lambda, \mu \in H^0(S, {\o}_S)^*$, $a_i, b_j \in H^0(S, {\o}_S)$,
$v(u) = n$, $v(w) = m$, and $a_{-i}, b_{-j}$ are nilpotent
elements of $H^0(S, {\o}_S)$, one has that
$$(u,w)_{p} = (-1)^{n\cdot m}
 \left ( \frac {{\lambda}^m
{\prod}_{i=1}^{\infty}{\prod}_{j=1}^{h}\big (1 -
a_i^{j/(i,j)}b_{-j}^{i/(i,j)})^{(i,j)}} {{\mu}^n
{\prod}_{j=1}^{\infty}{\prod}_{i=1}^{l}\big (1 -
b_j^{i/(i,j)}a_{-i}^{j/(i,j)})^{(i,j)}}\right )\, ,$$ \noindent
where, finitely, many of the terms appearing in the products
differ from 1 in ${\Bbb G}_m^{\bullet} (S) = H^0(S, {\o}_S)^*$. This
symbol is a generalization of the Contou-Carr\`{e}re symbol
associated with the rational point $p$.
\end{cor}

\begin{cor}\label{21} If we add the condition char($k$) = 0 to the hypothesis of Corollary \ref{c:2003},
we have that
 $$(u,w)_{p} = (-1)^{n\cdot m}
 \left ( \frac {{\lambda}^m
\cdot \exp (\underset {i>0} \sum (\mod {{\delta}_i (u)\cdot
{\delta}_{-i} (w)}/{i}))} {{\mu}^n \cdot \exp (\underset {i>0}
\sum (\mod {{\delta}_{-i} (u)\cdot {\delta}_{i} (w)}/{i}))}\right
)\, ,$$ where ${\delta}_s (f) = \res (z^s \cdot \frac {df}{f})$.
This expression was offered by C. Contou-Carr\`{e}re in \cite{Co},
Corollary 4.5.
\end{cor}

\begin{rem} By considering $\sp k$-valued points in the above morphism of
schemes, if $f,g \in {\Sigma}_C^*$ one has that $$(f,g)_p =
(-1)^{v_p(f)\cdot v_p(g)} \frac {f^{v_p(g)}}{g^{v_p(f)}} (p)\, ,
$$\noindent which is the expression of the multiplicative local symbol
defined by J.-P. Serre \cite{Se}. We should note that the
definition offered is ``local"  because we do not need to use the
reciprocity law $\underset {p\in C} \prod (f,g)_p = 1$ to
determine the uniqueness of the symbol. Indeed, we only need a
closed point $p\in C$ of a complete, irreducible and non-singular
curve to define it. Moreover, the conditions that appear in the
definition are natural in the theory of local symbols.
\end{rem}

\subsection{Generalization of
the Contou-Carr\`{e}re symbol associated with a closed point}

Let us now consider a complete, irreducible and non-singular curve
$C$ over a perfect field $k$ and with a rational point, and let
$p\in C$ an arbitrary closed point. Analogously to the previous
section, we set
$$K^{\Theta}_{\{p\}} = \underset {n} \varprojlim
K[(\tilde {{\mathcal L}_{\Theta}})^{n)}]\, .$$
        If $k(p)$ is the residue class field of the closed point $p$,
it follows from Theorem \ref{t:323} and Corollary \ref{cor:1234}
that $K^{\Theta}_{\{p\}}\simeq {\mathbb Z}_* \times
({\Gamma}_+)_{k(p)}$, with
$$({\Gamma}_+)_{k(p)}^{\bullet} (S) = {\Gamma}_{+}^{\bullet}
(S\underset {\sp k} \times \sp k(p))$$\noindent for each
$k$-scheme $S$.

    Let $({\Gamma}_{-})_{k(p)}$ and $({\mathbb G}_m)_{k(p)}$ be the
group schemes constructed in \cite{Pa3}, whose functors of points
are
$$\aligned ({\Gamma}_{-})_{k(p)}^{\bullet}(S) &= {\Gamma}_{-}^{\bullet}
(S\underset {\sp k} \times \sp k(p)) \\
({\mathbb G}_m)_{k(p)}^{\bullet}(S) &= {\mathbb G}_{m}^{\bullet}
(S\underset {\sp k} \times \sp k(p)) \endaligned$$ \noindent for
each $k$-scheme $S$. We can thus consider the group scheme ${\Bar
K}^{\Theta}_{\{p\}} = ({\mathbb G}_m)_{k(p)} \times
({\Gamma}_{-})_{k(p)}$ and set ${\tilde {\Gamma}}_p =
K^{\Theta}_{\{p\}} \times {\Bar K}^{\Theta}_{\{p\}}$, which is a
locally connected group scheme. Also, for each $k$-scheme $S$ one
has that
$${\tilde {\Gamma}}_p^{\bullet} (S) \simeq (H^{0}(S,{\o}_S)\underset {k}
\otimes k(p)) ((z))^{*}\, ,$$\noindent and therefore ${\tilde
{\Gamma}}_p^{\bullet} (\sp k) \simeq k(p)((z))^{*} \simeq (\Hat
{{\o}_p})^{*}_{(0)}$.

        Furthermore, from the natural immersion of functors of
groups
$${\mathbb Z}_*^{\bullet} (S) \hookrightarrow ({\mathbb
Z}_*)_{k(p)}^{\bullet} (S) =  {\mathbb Z}_*^{\bullet} (S\underset
{\sp k} \times \sp k(p))\ ,$$\noindent and the duality referred to
in Subsection \ref{ss:41}, one has defined a morphism of schemes:
$${\chi}_p \colon K^{\Theta}_{\{p\}} \times {\Bar
K}^{\Theta}_{\{p\}} \longrightarrow ({\mathbb G}_m)_{k(p)}.$$

    If we now denote by $N_{k(p)/k}\colon ({\mathbb
G}_m)_{k(p)}\to  {\mathbb G}_m$ the morphism of group schemes
induced by the norm of the extension of fields $k\hookrightarrow
k(p)$, we have the morphism:
$${\bar \chi} \colon K^{\Theta}_{\{p\}}
\times {\Bar K}^{\Theta}_{\{p\}} \longrightarrow {\mathbb G}_m\, ,
$$\noindent defined by
$${\bar \chi} (f,g)
= N_{k(p)/k}({\chi}_p(f,g)).$$

    Using ${\Bar {\mathcal H}} ({\tilde {\Gamma}}_p)$ to denote
the group scheme ${\mathbb G}_m \times K^{\Theta}_{\{p\}} \times
{\Bar K}^{\Theta}_{\{p\}}$ characterized by the group law
$$(\alpha,f,g)\cdot ({\alpha}',f',g')= ({\alpha}\cdot
{\alpha}'\cdot {\bar \chi}(f',g),f\cdot f',g\cdot g')$$\noindent
for $S$-valued points, $S$ being a $k$-scheme, keeping the
notations of the above section and bearing in mind that ${\tilde
{\Gamma}}_p = {\tilde {{\Gamma}_{k(p)}}}$, with ${\tilde
{{\Gamma}_{k(p)}}}$ the group scheme constructed in \cite{Pa3},
Section 2.B, and with a equivalent definition of $v$ to the
offered above, one has that:

\begin{thm}\label{20} There exists a unique element $(\quad,
\quad)_{p}$ in the cohomology class $[e_{{\Bar {\mathcal H}}
 ({\tilde {\Gamma}}_p)}] \in H_{reg}^2
({\tilde {\Gamma}}_p,{\mathbb G}_m)$ satisfying the conditions:

\begin{itemize}
\item $(f,g\cdot g')_p = (f,g)_{p} \cdot (f,g')_{p}$;
\item $(f,g)_{p} = e_{{\Bar {\mathcal H}}
({{{\tilde {\Gamma}}_{p}}})} (f,g)$ if  $v(f)= 0$;
\item $(f,-f)_{p} = 1$
\end{itemize}
\noindent for $f,g,g' \in {{\Gamma}_p}^{\bullet} (S)$, with $S$ a
connected $k$-scheme. This element is a generalization of the
Contou-Carr\`{e}re symbol associated with the closed point $p$.
\end{thm}

\begin{cor}\label{c:nor} If $S$ is a connected
$k$-scheme and $u,w \in {\tilde {\Gamma}}_p^{\bullet} (S)$ with
$$\begin{aligned} \quad &u = \lambda z^n{\prod}_{i=1}^{l}
(1-a_{-i}z^{-i}) {\prod}_{i=1}^{\infty} (1-a_{i}z^{i}) \\
&w = \mu z^m {\prod}_{j=1}^{h} (1-b_{-j}z^{-j})
{\prod}_{j=1}^{\infty} (1-b_{j}z^{j})\, ,
\end{aligned}$$\noindent where $\lambda, \mu \in [H^0(S, {\o}_S)\otimes_k k(p)]^*$,
$a_i, b_j \in H^0(S, {\o}_S)\otimes_k k(p)$,
$v(u) = n$, $v(w) = m$, and $a_{-i}, b_{-j}$ are nilpotent
elements of $H^0(S, {\o}_S)\otimes_k k(p)$, one has that
$$(u,w)_{p} = (-1)^{n\cdot m\cdot deg(p)} N_{k(p)/k}
 \left ( \frac {{\lambda}^m
{\prod}_{i=1}^{\infty}{\prod}_{j=1}^{h}\big (1 -
a_i^{j/(i,j)}b_{-j}^{i/(i,j)})^{(i,j)}} {{\mu}^n
{\prod}_{j=1}^{\infty}{\prod}_{i=1}^{l}\big (1 -
b_j^{i/(i,j)}a_{-i}^{j/(i,j)})^{(i,j)}}\right )\, ,$$ \noindent
where, finitely, many of the terms appearing in the products
differ from 1 in ${\Bbb G}_m^{\bullet} (S) = H^0(S, {\o}_S)^*$.
\end{cor}

\begin{cor}\label{3421} If we add the condition char($k$) = 0 to the hypothesis of the previous corollary,
we have that
 $$(u,w)_{p} = (-1)^{n\cdot m \cdot deg(p)} N_{k(p)/k}
 \left ( \frac {{\lambda}^m
\cdot \exp (\underset {i>0} \sum (\mod {{\delta}_i (u)\cdot
{\delta}_{-i} (w)}/{i}))} {{\mu}^n \cdot \exp (\underset {i>0}
\sum (\mod {{\delta}_{-i} (u)\cdot {\delta}_{i} (w)}/{i}))}\right
)\, ,$$ where ${\delta}_s (f) = \res (z^s \cdot \frac {df}{f})$.
\end{cor}

\subsection{Generalized Contou-Carr\`{e}re symbol associated with a morphism}

    Keeping the notations of the previous subsection, let
$\varphi\colon {\mathbb G}_m \to G$ be a morphism of group
schemes, $G$ being an algebraic group scheme.

    Let us now consider the group scheme ${\Bar {\mathcal H}_{\varphi}} ({\tilde {\Gamma}}_p) =
G\times K^{\Theta}_{\{p\}} \times {\Bar K}^{\Theta}_{\{p\}}$,
characterized by the group law
$$(a,f,g)\cdot (a',f',g')= (a\cdot
a'\cdot {\bar \chi}_{\varphi}(f',g),f\cdot f',g\cdot
g')$$\noindent for $S$-valued points on a connected $k$-scheme,
and with ${\bar \chi}_{\varphi} = \varphi \circ {\bar \chi}$.

    If we denote by $e_{{\Bar {\mathcal H}_{\varphi}}
 ({\tilde {\Gamma}}_p)}$ the commutator of the induced
  extension of group schemes
 $$0 \to G \to {\Bar {\mathcal H}_{\varphi}} ({\tilde
 {\Gamma}}_p)\to {\tilde
 {\Gamma}}_p\to 0\, ,$$\noindent similar to the above cases, one has that:

\begin{thm}\label{t:morph} There exists a unique element $(\quad,
\quad)_{p}^{\varphi}$ in the cohomology class $[e_{{\Bar {\mathcal
H}_{\varphi}}
 ({\tilde {\Gamma}}_p)}] \in H_{reg}^2
({\tilde {\Gamma}}_p,G)$ satisfying the conditions:

\begin{itemize}
\item $(f,g\cdot g')_p^{\varphi} = (f,g)_{p}^{\varphi} \cdot (f,g')_{p}^{\varphi}$;
\item $(f,g)_{p}^{\varphi} = e_{{\Bar {\mathcal H}_{\varphi}}
({{{\tilde {\Gamma}}_{p}}})} (f,g)$ if  $v(f)= 0$;
\item $(f,-f)_{p}^{\varphi} = 1$
\end{itemize}
\noindent for $f,g,g' \in {{\Gamma}_p}^{\bullet} (S)$, with $S$ a
connected $k$-scheme. This element is a generalization of the
Contou-Carr\`{e}re symbol associated with the morphism
$\varphi\colon {\mathbb G}_m \to G$.
\end{thm}

\begin{cor}\label{c:nor2} If $S$ is a connected
$k$-scheme and $u,w \in {\tilde {\Gamma}}_p^{\bullet} (S)$ with
$$\begin{aligned} \quad &u = \lambda z^n{\prod}_{i=1}^{l}
(1-a_{-i}z^{-i}) {\prod}_{i=1}^{\infty} (1-a_{i}z^{i}) \\
&w = \mu z^m {\prod}_{j=1}^{h} (1-b_{-j}z^{-j})
{\prod}_{j=1}^{\infty} (1-b_{j}z^{j})\, ,
\end{aligned}$$\noindent where $\lambda, \mu \in [H^0(S, {\o}_S)\otimes_k k(p)]^*$,
$a_i, b_j \in H^0(S, {\o}_S)\otimes_k k(p)$, $v(u) = n$, $v(w) = m$,
and $a_{-i}, b_{-j}$ are nilpotent elements of $H^0(S,
{\o}_S)\otimes_k k(p)$, one has that
$$(u,w)_{p}^{\varphi} = (g_{-1})^{n\cdot m\cdot deg(p)} \varphi \left (N_{k(p)/k}
 \left ( \frac {{\lambda}^m
{\prod}_{i=1}^{\infty}{\prod}_{j=1}^{h}\big (1 -
a_i^{j/(i,j)}b_{-j}^{i/(i,j)})^{(i,j)}} {{\mu}^n
{\prod}_{j=1}^{\infty}{\prod}_{i=1}^{l}\big (1 -
b_j^{i/(i,j)}a_{-i}^{j/(i,j)})^{(i,j)}}\right )\right ) ,$$
\noindent where, finitely, many of the terms appearing in the
products differ from 1 in $G^{\bullet} (S)$, and $g_{-1} = \varphi
(-1)$.
\end{cor}

\begin{cor}[Reciprocity Law]\label{c:law} Let $C$ be a complete, irreducible,
non-singular, and with a rational point curve over a perfect field
$k$. For each morphism of group schemes $\varphi\colon {\mathbb
G}_m \to G$, and each artinian local finite $k$-algebra $B$, one
has that
$$\underset {p\in C} \prod (f,g)_{p}^{\varphi} = 1 \text { for all
} f, g \in ({\Sigma}_C\otimes_k B)^*\, ,$$\noindent where,
finitely, many of the terms appearing in the products differ from
1 in $G^{\bullet} (\sp B)$.
\end{cor}

\begin{proof} From the natural immersion $({\Sigma}_C\otimes_k B)^* \hookrightarrow
{\tilde {\Gamma}}_p^{\bullet} (\sp B)$, and following the method
of \cite{AP}, the claim is a direct consequence of  the statements
proved in \cite{Pa3}.
\end{proof}

    Henceforth, $C$ will be an irreducible, complete,
non-singular and with a rational point curve over a finite field
$k$ that contains the $m^{th}$-roots of unity, with $\# k = q$.

    Let us now consider $\phi_N$, the character associated with the integer number $N\in
{\mathbb Z}$,
$$\begin{aligned} \phi_N\colon {\mathbb G}_m
&\longrightarrow {\mathbb G}_m \\ \lambda &\longmapsto {\lambda}^N
\end{aligned}$$\noindent where $\lambda \in {\mathbb G}_m^{\bullet}
(S)$, with $S$ a connected $k$-scheme.

\begin{lem}\label{l:art} For each closed point $p\in C$ and each
artinian local finite $k$-algebra $B$, there exists a positive
integer number $N$ such that the morphism of group schemes
$$(\quad,
\quad)_{p}^{\phi_N}\colon {{\Gamma}_p} \times {{\Gamma}_p}
\longrightarrow {\mathbb G}_m$$\noindent satisfies a factorization
    $$\xymatrix{{\tilde {\Gamma}}_p^{\bullet}(\sp B)\times
    {\tilde {\Gamma}}_p^{\bullet}(\sp B)\ar[rr]^-{(\quad,\quad)_{p}^{\phi_N}}
\ar[dr]_-{\tilde {(\quad, \quad)}_{p,B}^{\phi_N}} & & {\mathbb
G}_m^{\bullet} (\sp B) \\ & {\mathbb \mu}_m^{\bullet} (\sp B)
\ar[ru]_{i_{B}}& \, ,}$$\noindent where ${\mathbb \mu}_m$ is the
$k$-scheme of $m^{th}$-roots of unity and $i_{B}$ is the natural
immersion.
\end{lem}

\begin{proof} If $\alpha = \sharp B^*$, since ${\mathbb \mu}_m^{\bullet} (\sp B) =
\{b\in B^* \text{ such that } b^m = 1\}$, and since $m$ divides
$\alpha$ because $k^*$ is a subgroup of $B^*$, the statement
follows from considering the integer number $N = \frac
{\alpha}{m}$.
\end{proof}

    One has that:

\begin{prop}\label{p:charact}${\tilde {(\quad,
\quad)}_{p,B}^{\phi_N}}$ is the unique element in the cohomology
class $[e_{{\Bar {\mathcal H}}_{\phi_N}({\tilde
{\Gamma}}_p^{\bullet}(\sp B))}] \in H^2 ({\tilde
{\Gamma}}_p^{\bullet}(\sp B),{\mathbb \mu}_m^{\bullet} (\sp B))$
satisfying the conditions:

\begin{itemize}
\item $\tilde {(f,g\cdot g')}_{p, B}^{\phi_N} = \tilde {(f,g)}_{p,
B}^{\phi_N} \cdot \tilde {(f,g')}_{p, B}^{\phi_N}$; \item $\tilde
{(f,g)}_{p, B}^{\phi_N} = e_{{\Bar {\mathcal H}}_{\phi_N}({\tilde
{\Gamma}}_p^{\bullet}(\sp B))} (f,g)$ if  $v(f)= 0$; \item $\tilde
{(f,-f)}_{p, B}^{\phi_N} = 1$
\end{itemize}
\noindent for $f,g,g' \in {\tilde {\Gamma}}_p^{\bullet}(\sp B)$.
\end{prop}

\begin{proof} The claim follows from the definition of ${\tilde {(\quad,
\quad)}_{p, B}^{\phi_N}}$ and the above factorization.
\end{proof}

\begin{rem}  If $B = k$, and we consider the
immersion of groups $${\Sigma}_C^* \hookrightarrow {\tilde
{\Gamma}}_p^{\bullet}(\sp k)\, ,$$\noindent the symbol ${\tilde
{(\quad, \quad)}_{p, k}^{\phi_{\frac {(q-1)}{m}}}}$ is the Hilbert
norm residue symbol for the closed point $p\in C$. Hence the
morphism of schemes $(\quad, \quad)_{p}^{\phi_N}$ is a
generalization of the construction of the Hilbert norm residue
symbol made in \cite{Pa0}.
\end{rem}

\begin{cor}\label{c:rec-law} If $f, g \in ({\Sigma}_C\otimes_k
B)^*$, one has that $$\underset {p\in C} \prod  {\tilde
{(f,g)}_{p, B}^{\phi_N}} = 1\, ,$$\noindent where, finitely, many
of the terms appearing in the products differ from 1 in ${\mathbb
\mu}_m^{\bullet} (\sp B)$.
\end{cor}

\begin{proof} This reciprocity law is a direct consequence of the
factorization referred to in Lemma \ref{l:art} and the reciprocity
law of Corollary \ref{c:law}.
\end{proof}

\begin{rem} If $B = k$, the formula $$\underset {p\in C} \prod  {\tilde
{(f,g)}_{p, k}^{\phi_{\frac {(q-1)}{m}}}} = 1$$\noindent is the
reciprocity law of the Hilbert norm residue symbol.
\end{rem}

\bigskip {\centerline {\bf ACKNOWLEDGMENT}} The main results of the
present paper are part of the author's Ph.D. thesis. I would like
to express my gratitude to my thesis advisor, Prof. Dr. Jos\'e
Mar\'{\i}a Mu\~noz Porras, for his constant support and help with
this work.

\end{document}